\documentclass[10pt, a4paper]{article}

\usepackage{amsmath,  amsfonts}

\usepackage[T1]{fontenc}

\usepackage{cite}

\newtheorem{theorem}{\quad Theorem}[section]

\newcommand{\be} {\begin{equation}}
\newcommand{\ee} {\end{equation}}

\title {An interior boundedness result for an elliptic equation.}

\date{}

 \author{Samy Skander Bahoura\footnote {e-mails: samybahoura@yahoo.fr, samybahoura@gmail.com} \\ 
 {\small Equipe d'Analyse Complexe et G\'eom\'etrie.}\\  
  {\small Universit\'e Pierre et Marie Curie, 75005 Paris, France.}} 

\begin{document}

\maketitle
 
\begin{abstract}
We derive a local unfiorm boundedness result for an equation with weight having interior singularity.

\end{abstract}



\section{Introduction and Main Results} 

We set $ \Delta = \partial_{11} + \partial_{22} $  on open set $ \Omega $ of $ {\mathbb R}^2 $ with a smooth boundary.

\bigskip

We consider the following equation:

$$ (P)   \left \{ \begin {split} 
      -\Delta u & = \dfrac{1}{-\log \dfrac{|x|}{2d}}V e^{u} \,\, &\text{in} \,\, & \Omega  \subset {\mathbb R}^2, \\
                  u & = 0  \,\,             & \text{in} \,\,    &\partial \Omega.              
\end {split}\right.
$$

Here:

$$ 0\leq V\leq b,\,\, \int_{\Omega} \dfrac{1}{-\log \dfrac{|x|}{2d}} e^u dx \leq C, \,\, u \in W_0^{1,1}(\Omega), $$

and,

$$ d=diam(\Omega), \,\, 0 \in \Omega $$

Equations of the previous type were studied by many authors, with or without  the boundary condition, also for Riemannian surfaces,  see ~\cite{1,2,3,4,5,6,7,8,9,10,11,12,13,14,15,16,17,18,19},  where one can find some existence and compactness results.

Among other results, we  can see in \cite{11} the following important Theorem

\smallskip

{\bf Theorem A}{\it (Brezis-Merle \cite{11})}.{\it If $ (u_i) $ is a sequence of solutions of problem $ (P) $ with $ (V_i) $ satisfying $ 0 < a  \leq V_i  \leq b < + \infty $ and without the term $ \dfrac{1}{-\log \dfrac{|x|}{2d}} $, then, for any
compact subset $ K $ of $  \Omega $, it holds:

$$ \sup_K u_i \leq c, $$ 

with c depending on $ a, b, K, \Omega $}

One can find in \cite{11} an interior estimate if we assume $ a=0 $, but we need an assumption on the integral of $ e^{u_i} $, namely, we have:

\smallskip

{\bf Theorem B}{\it (Brezis-Merle \cite{11})}.{\it For $ (u_i)_i $ and $ (V_i)_i $ two sequences of functions relative to the problem $ (P) $ without the term $ \dfrac{1}{-\log \dfrac{|x|}{2d}} $ and  with,
$$ 0 \leq V_i \leq b < + \infty \,\, {\rm and} \,\, \int_{\Omega} e^{u_i} dy  \leq C, $$
then for all compact set $ K $ of $ \Omega $ it holds;
$$ \sup_K u_i \leq c, $$
with $ c $ depending on $ b, C, K $ and $ \Omega $.}

\smallskip

If we assume $ V $ with more regularity, we can have another type of estimates, a $ \sup + \inf $ type inequalities. It was proved by Shafrir see \cite{18}, that, if $ (u_i)_i $ is a sequence of functions solutions of the previous equation without assumption on the boundary  with $ V_i $ satisfying $ 0 < a \leq V_i \leq b < + \infty $, then we have a $ \sup + \inf $ inequality.

\smallskip

Here, we have:

\smallskip

\begin{theorem}For sequences $ (u_i)_i $ and $ (V_i)_i $ of the Problem $ (P) $, for all compact subsets $ K $ of $ \Omega $ we have:

$$  || u_i||_{L^{\infty}(K)} \leq c(b, C, K, \Omega), $$

\end{theorem}

{\bf Remark:} Remark that we have a $ C^0 $ weight $ \dfrac{1}{-\log \dfrac{|x|}{2d}} $, the solutions are not $ C^2 $, but if we add some assumptions on $ V_i $ one can consider $ C^2 $ solutions and $ C^2 $ convergence of sequences.

\smallskip

On can have the regularity $ C^2 $ of the solutions and the $ C^2 $ convergence of the solutions if we suppose for example $ V_i\in C^{0,\epsilon}, \epsilon >0 $ and for the convergence $ V_i\to V $ in the space $ C^{0,\epsilon}, \epsilon >0 $. Indeed, one can reduce the problem to regularity and convergence of  the Newtonian potential of a radial distribution $ f(x)=f(|x|)=\dfrac{V_i(0)e^{u_i(0)}}{-\log \left (\dfrac{|x|}{2d}\right )} \eta (|x|) $, with $ \eta $ a cutoff function ($ \eta \equiv 1 $ in a neighborhood of $ 0 $ with compact support and radial), see for example the book of Dautray-Lions, chapter 2, Laplace operator.

\smallskip

By a duality theorem one can prove that (see \cite{12}):

$$ ||\nabla u_i||_q \leq C_q, \,\, \forall \, 1 \leq q <2. $$

If we add the assumption that

$$ ||\nabla V_i||_{\infty} \leq A, $$
 
then by a result of Chen-Li of "moving-plane" we have a compactness of $ (u_i)_i $ near the boundary, see  \cite{13}.

\smallskip

We ask the following question about inequality of type $ \sup + \inf $, as in the work of Tarantello, see \cite{19} and Bartolucci-Tarantello, see \cite{8}:

{\bf Problems}. 1) Consider the Problem $ (P) $ without the boundary condition (without Dirichlet condition) and assume that:

$$ 0 < a \leq V \leq b <+\infty, $$

Does exists constants $ C_1=C_1(a, b, K, \Omega), C_2=C_2(a, b, K, \Omega) $ such that:

$$ \sup_K u + C_1 \inf_{\Omega} u \leq C_2, $$

for all solution $ u $ of $ (P) $ ?

\smallskip

2) If we add the condition $ ||\nabla V ||_{\infty} \leq A $, can we have a sharp inequality:

$$ \sup_K u + \inf_{\Omega} u \leq c(a,b, A, K, \Omega) ? $$

\section{Proof of the Theorem} 

We have:

$$ u_i \in W_0^{1,1}(\Omega),\,\, {\rm and} \,\,\dfrac{1}{-\log \dfrac{|x|}{2d}} e^{u_i}\in L^1(\Omega). $$

Thus, by corollary 1 of Brezis and Merle we have:

$$ e^{u_i} \in L^k(\Omega), \forall \, k >2 . $$

Using the elliptic estimates and the Sobolev embedding, we have:

$$ u_i \in W^{2,k}(\Omega) \cap C^{1,\epsilon}(\bar \Omega). $$

By the maximum principle $ u_i \geq 0 $.

\smallskip

Also, by a duality theorem or a result of Brezis-Strauss, we have:

$$ ||\nabla u_i||_q \leq C_q, \,\, 1 \leq q < 2. $$

Since,

$$ \int_{\Omega} \dfrac{1}{-\log \dfrac{|x|}{2d}} V_i e^{u_i} dx \leq C, $$

We have a convergence to a nonegative measure $ \mu $:

$$ \int_{\Omega} \dfrac{1}{-\log \dfrac{|x|}{2d}} V_i e^{u_i} \phi dx \to \int_{\Omega} \phi d\mu, \,\, \forall \,\, \phi \in C_c(\Omega). $$

We set $ S $ the following set:

\smallskip

$ S= \{ x \in  \Omega , \exists \, (x_i) \in \Omega, \, x_i  \to x, \, u_i (x_i)  \to + \infty \} $.

\smallskip

We say that $ x_0 $ is a regular point of $ \mu $ if there function $ \psi \in C_c(\Omega) $, $ 0 \leq \psi \leq 1 $, with $ \psi =1 $ in a neighborhood of $ x_0 $ such that:

\be \int \psi d\mu < 4 \pi. \label{(1)}  \ee

We can deduce that a point $ x_0 $ is non-regular if and only if $ \mu({x_0}) \geq 4\pi . $

\smallskip

A consequence of this fact is that if $ x_0 $ is a regular point then:

\be \exists \,\, R_0 > 0 \,\, such \,\, that \,\, one \, can\,\, bound  \,\, (u_i)=(u_i^+) \,\, in  \,\, L^{\infty}(B_{R_0}(x_0)).\label{(2)} \ee 

We deduce $ (\ref{(2)})$ from corollary 4 of Brezis-Merle paper, because we have by the Gagliardo-Nirenberg-Sobolev inequality:

$$ || u_i^+||_1=||u_i||_1\leq c_q||u_i||_{q^*}\leq C_q'||\nabla u_i||_q \leq C_q, \,\, 1 \leq q < 2. $$

We denote by $ \Sigma $ the set of non-regular points.

\underbar { Step 1}: S = $ \Sigma $.

We have $  S \subset \Sigma $. Let's consider $ x_0 \in \Sigma $. Then we have:

\be   \forall \,\,\,  R >0,\,\,\,  \lim ||u_i^+ ||_{ L^{\infty }(B_R(x_0))} = + \infty. \label{(3)} \ee

Suppose contrary that:

$$  ||u_i^+ ||_{ L^{\infty }(B_{ R_0}(x_0))} \leq C. $$

Then:

$$  ||e^{u_{i_k}} ||_{ L^{\infty }(B_{ R_0}(x_0))} \leq C, \,\,\, {\rm and  }$$

$$ \int_{B_{ R}(x_0)} \dfrac{1}{-\log \dfrac{|x|}{2d}}    V_{i_k}e^{u_{i_k}} = o(1). $$

For $ R $ small enough, which imply $ (\ref{(1)}) $ for a function $ \psi $  and $ x_0 $ will be  regular, contradiction. Then we have $ (\ref{(3)}) $. We choose  $ R_ 0 > 0  $ small such that $ B_{ R_0}(x_0) $ contain only $ x_0 $ as non -regular point. $ \Sigma $. Let's $ x_i \in B_R(x_0)  $ scuh that:

$$ u_i^+(x_i) = \max_{B_R(x_0)} u_i^+ \to + \infty. $$

We have $ x_i \to x_0 $. Else, there exists $ x_{i_k} \to \bar x \not = x_0 $ and $ \bar x  \not \in \Sigma $, i.e. $ \bar x $ is a regular point. It is impossible because we would have $ (\ref{(2)}) $. 

\smallskip

Since the measure is finite, if there are blow-up points, or non-regular points, $ S =\Sigma $ is finite.

\smallskip

\underbar { Step 2}: $ \Sigma =\{ \emptyset \} $.

\smallskip

Now: suppose contrary that there exists a non-regular point $ x_0 $. We choose a radius $ R >0 $ such that $ B_R(x_0) $ contain only $ x_0 $ as non-regular point. Thus outside $ \Sigma $ we have local unfirorm boundedness of $ u_i $, also in $ C^1 $ norm. Also, we have weak *-convergence of $ V_i $ to $ V \geq 0 $ with $ V \leq b $.

\smallskip

Let's consider (by a variational method):

$$ z_i \in W_0^{1,2}(B_{R}(x_0)), $$

$$ -\Delta z_{i} = f_{i}=\dfrac{1}{-\log \dfrac{|x|}{2d}} V_i e^{u_i}\,\, {\rm in } \,\,\, B_R(x_0),\,\,\,{\rm et} \,\,\, z_{i} = 0 \,\,\,{\rm on } \,\,\, \partial  B_R(x_0). $$

By a duality theorem:

$$ z_i \in W_0^{1,q}(B_{R}),\,\, ||\nabla z_i||_q \leq C_q. $$

By the maximum principle, $ u_{i} \geq z_{i} $ in $ B_R(x_0) $.

\be \int \dfrac{1}{-\log \dfrac{|x|}{2d}} e^{z_i} \leq
\int \dfrac{1}{-\log \dfrac{|x|}{2d}} e^{u_{i}}  \leq C. \label{(*)} \ee

On the other hand, $ z_{i} \to z $ a.e.  (uniformly on compact sets of $ B_R(x_0) -  \{ x_0\} $) with $ z $ solution of :

$$ -\Delta z = \mu \,\,\, {\rm in } \,\,\, B_R(x_0),\,\,\,{\rm et} \,\,\, z = 0 \,\,\,{\rm on } \,\,\,  \partial B_R(x_0). $$

Also, we have up to a subsequence, $ z_i \to z $ in $ W_0^{1,q}(B_R(x_0)), 1\leq q < 2 $ weakly, and thus $ z \in W_0^{1,q}(B_R(x_0)) $.

\smallskip

Then by Fatou lemma:

\be \int \dfrac{1}{-\log \dfrac{|x|}{2d}} e^{z} \leq C. \label{(**)} \ee

As $ x_0 \in S $  is not regular point we have $ \mu(\{ x_0\}) \geq 4 \pi $, which imply that,  $ \mu \geq 4 \pi \delta_{x_0} $ and by the maximum principle in $ W^{1,1}_0(B_R(x_0)) $ (obtainded by Kato's inequality)

$$ z(x) \geq  2 \log \dfrac{1}{|x - x_0|} + O(1) \, \, {\rm if } \, \, x  \to x_0. $$

Because,

$$ z_1\equiv2 \log \dfrac{1}{|x - x_0|}+ 2\log R \in W_0^{1,s}(B_{R}(x_0)),\,\, 1\leq s < 2. $$

Thus,

$$ \dfrac{1}{-\log \dfrac{|x|}{2d}} e^{z} \geq  \dfrac{C}{-|x-x_0|^2\log \dfrac{|x|}{2d}}, \, \, C > 0. $$

Both in the cases $ x_0=0 $ and $ x_0 \not = 0 $ we have:

$$  \int_{B_R(x_0)}\dfrac{1}{-\log \dfrac{|x|}{2d}} e^z = \infty. $$

But, by $ (\ref{(**)}) $:

$$ \int \dfrac{1}{-\log \dfrac{|x|}{2d}} e^z   \leq C. $$

which a contradiction.

\end{document}